\documentclass[final,10pt]{article}
\usepackage{verbatim}
\usepackage{latexsym}
\usepackage{amsmath,amsthm,amsfonts}
\usepackage{eucal}
\usepackage[notref,notcite]{showkeys}

\usepackage{mathtext}
\usepackage[cp1251]{inputenc}
\usepackage[T2A]{fontenc}
\usepackage[dvips]{graphicx}
\usepackage{amsmath}
\usepackage{amssymb}
\usepackage{amsxtra}
\usepackage{latexsym}
\usepackage{ifthen}

\newtheorem{theorem}{Theorem}
\newtheorem{definition}{Definition}[section]
\newtheorem{lemma}{Lemma}[section]

\newtheorem{remark}{Remark}[section]

\numberwithin{equation}{section}

\DeclareMathOperator{\supp}{supp}

\begin{document}

\title {\bf  Finite speed of propagations of the electromagnetic field
in nonlinear isotropic dispersive mediums}

\author{{\normalsize\bf Yuliya V. Namlyeyeva, Roman M. Taranets}
\footnote{Research is partially supported by the INTAS project
Ref. No: 05-1000008-7921}
\smallskip}

\def\lhead{R.M. Taranets}

\def\rhead{Maxwell equations}

\date{\today}

\maketitle

\setcounter{section}{0}

\begin{abstract}
We propose some modification of Maxwell's equations describing
mediums which electric and magnetic properties are changed
essentially after interaction with outer electromagnetic field. We
show for such mediums that electromagnetic waves have finite speed
of propagations property for some time depending on initial energy
of electromagnetic field and nonlinear parameters of the problem
which are responsible for properties of medium.
\end{abstract}

\textbf{2000 MSC:} {35Q60, 35B40, 78A25}

\textbf{keywords:} {Maxwell's equations, nonlinear dispersive
medium, finite speed of propagations, asymptotic behavior}

\section{Introduction}

We consider classical Maxwell system (see \cite{LL}):
\begin{gather} \qquad   \tfrac{1}{c}\textbf{D}_t + \tfrac{4\pi}{c}\textbf{J}
= \text{curl}\,\textbf{H}, \hfill \label{01}\\
\tfrac{1}{c}\textbf{B}_t  + \text{curl}\,\textbf{E} = 0, \hfill \label{02}\\
\text{div} \textbf{D} = 4 \pi \rho,\ \text{div} \textbf{B} = 0,
\hfill \label{03}
\end{gather}
\vspace{-2.2cm}
$$
(M_0)\qquad \left\{ \hspace{+3cm}\phantom {\begin{gathered} \qquad
\tfrac{1}{c}\textbf{D}_t + \tfrac{4\pi}{c}\textbf{J}
= \text{curl}\,\textbf{H}, \hfill \label{01}\\
\tfrac{1}{c}\textbf{B}_t  + \text{curl}\,\textbf{E} = 0, \hfill \label{02}\\
\text{div} \textbf{D} = 4 \pi \rho,\ \text{div} \textbf{B} = 0,
\label{03}, \hfill
\end{gathered}} \right.
$$
where $\textbf{E}$ and $\textbf{H}$ are electric and magnetic
fields; $\textbf{D}$ and $\textbf{B}$ are electric and magnetic
inductions; $\rho$ is charge density and $c$ is velocity of light.
The current density $\textbf{J}$  satisfies by Ohm's law:
\begin{equation}\label{Ohm}
\textbf{J} = \sigma \textbf{E},
\end{equation}
where $\sigma$ is electric conductivity.

We consider isotropic mediums in which permittivity $\varepsilon =
\varepsilon(x,t)$ and magnetic $\mu = \mu(x,t)$ conductivity are
functions of space and time. In this situation, state equations
have the following simple form (see \cite{LL}):
\begin{equation}\label{se}
\textbf{D} = \varepsilon \textbf{E},\  \textbf{B} = \mu
\textbf{H}.
\end{equation}
Substituting (\ref{se}) into equations (\ref{01}) and (\ref{02}),
we obtain equations for $\textbf{E}$ and $\textbf{H}$ for an
isotropic medium in the following dimensionless form:
\begin{gather} \qquad   \textbf{E}_t + a_1 \textbf{E}
- b_1 \text{curl}\,\textbf{H} = 0, \hfill \label{001}\\
\textbf{H}_t  + a_2 \textbf{H} + b_2 \text{curl}\,\textbf{E} = 0,
\hfill \label{002}
\end{gather}
\vspace{-1.6cm}
$$
(M_1)\qquad \left\{ \hspace{+3cm}\phantom {\begin{gathered} \qquad
\textbf{E}_t + a_1 \textbf{E}
- b_1 \text{curl}\,\textbf{H} = 0, \hfill \label{001}\\
\textbf{H}_t  + a_2 \textbf{H} + b_2 \text{curl}\,\textbf{E} = 0,
\hfill \label{002}
\end{gathered}} \right.
$$
where $a_i = a_i(x,t)$, $b_i = b_i(x,t)$, and
\begin{equation}\label{AB}
a_1 = \varepsilon^{-1} (\varepsilon_t + \sigma),\ a_2 = \mu^{-1}
\mu_t,\ b_1 = \varepsilon^{-1},\ b_2 = \mu^{-1}.
\end{equation}
Hyperbolic systems as ($M_1$) are well investigated (see, e.\,g.,
\cite{MS}).

In more general situation, electric and magnetic inductions depend
on electric and magnetic fields (see \cite{LL}), i.\,e.
\begin{equation}\label{se0}
\textbf{D} = \textbf{D}(\textbf{E},\textbf{H}), \  \textbf{B} =
\textbf{B}(\textbf{E},\textbf{H}).
\end{equation}
Below, we consider the simplest case of state equations
(\ref{se0}) when permittivity and magnetic conductivity  are some
functions of space and time depending on electric $E$ and magnetic
$H$ fields and its gradients, i.\,e. the relations (\ref{se}) with
$\varepsilon = \varepsilon(x,t, E, H, \nabla\,E, \nabla\,H)$ and
$\mu = \mu(x,t, E, H, \nabla\,E, \nabla\,H)$. In the case, we
arrive at the system ($M_1$), i.\,e. equations for $\textbf{E}$
and $\textbf{H}$ in an isotropic nonlinear medium, where $a_i=
a_i(x,t, E, H, \nabla\,E, \nabla\,H)$ and $b_i = b_i(x,t, E, H,
\nabla\,E, \nabla\,H)$ ($i=1,2$) satisfy relations (\ref{AB}). We
will study mediums in which nonlinear functions $a_i$ and $b_i$
satisfy the following conditions:
\begin{gather}
a_i(x,t, E, H, \nabla\,E, \nabla\,H) \geqslant d_1 w^{m-1}
|\nabla\,w|^p,
\ 0 < d_1 < \infty,\ i =1,2, \hfill \label{c1}\\
b_1(x,t, E, H,\nabla\,E, \nabla\,H) = b_2(x,t, E, H,\nabla\,E,
\nabla\,H),\label{c2}\hfill\\
|b_i(x,t, E, H,\nabla\,E, \nabla\,H)| \leqslant d_2 w^{n -1},\ 0 < d_2 < \infty,\ i =1,2, \label{c3}\hfill\\
|\nabla\,b_i(x,t, E, H, \nabla\,E, \nabla\,H)| \leqslant d_3 w^{n
-2}|\nabla\,w|, \ 0 < d_3 < \infty,\ i =1,2, \hfill \label{c4}
\end{gather}
where $w = w(x,t)= E^2 + H^2$ is dimensionless energy density
corresponding to isotropic mediums with constant permittivity and
magnetic conductivity;
\begin{equation}\label{mnp}
m \in \mathbb{R}^1, \ p > 0 \text{ and } n
> 0
\end{equation}
are parameters of medium. Conditions (\ref{c1})--(\ref{c4})
results in the following restrictions on $\varepsilon$ and $\mu$:
$$
\varepsilon = \mu \geqslant d_2^{-1}w^{1 - n},\ \varepsilon^{-1}
(\varepsilon_t + \sigma) \geqslant d_1 w^{m-1} |\nabla\,w|^p,\
\mu^{-1} \mu_t \geqslant d_1 w^{m-1} |\nabla\,w|^p,
$$
whence we deduce that
$$
\varepsilon = \mu \geqslant \max\{d_2^{-1}w^{1 - n},\
\varepsilon_{|_{t=0}} e^{ d_1 \int\limits_{0}^{t}{w^{m-1}
|\nabla\,w|^p\,d\tau}}\}.
$$
The equations like as $(M_1)$ describe mediums in which
permittivity and magnetic conductivity are some nonlinear
functions. The mediums have same structure to have to appear in
the simulation of various processes in laser optics and weakly
ionized plasma theory, where properties of medium are strongly
depend on energy density of electromagnetic field, for example,
ferroelectric, piezoelectric, multiferroic and etc.


In this paper, we study the propagation properties of solutions to
Cauchy problem for Maxwell's equations in the following
dimensionless form
\begin{gather} \qquad   \textbf{E}_t + a_1\,\textbf{E} - b_1\,\text{curl}\,\textbf{H}
= 0 \text{ in }Q_T, \hfill \label{1}\\
\textbf{H}_t  + a_2 \,\textbf{H} +
b_2\,\text{curl}\,\textbf{E} = 0 \text{ in }Q_T, \hfill \label{2}\\
\textbf{E}(0,x) = \textbf{E}_0 (x),\ \textbf{H}(0,x) =
\textbf{H}_0 (x), \hfill \label{4}
\end{gather}
\vspace{-2.14cm}
$$
(M)\qquad \left\{ \hspace{+3cm}\phantom {\begin{gathered}
\textbf{E}_t + a_1\,\textbf{E} - b_1\,\text{curl}\,\textbf{H}
= 0 \text{ in }Q_T, \hfill \label{1}\\
\textbf{H}_t  + a_2 \,\textbf{H} +
b_2\,\text{curl}\,\textbf{E} = 0 \text{ in }Q_T, \hfill \label{2}\\
\textbf{E}(0,x) = \textbf{E}_0 (x),\ \textbf{H}(0,x) =
\textbf{H}_0 (x), \hfill \label{4}
\end{gathered}} \right.
$$
where $Q_T = (0,T) \times \mathbb{R}^N$, $N = 2,3$, $0 < T
<\infty$, and the functions $a_i= a_i(x,t, E, H, \nabla\,E,
\nabla\,H)$, $b_i = b_i(x,t, E, H, \nabla\,E, \nabla\,H)$
($i=1,2$) satisfy conditions (\ref{c1})--(\ref{c4}). The unknown
functions are electric $\textbf{E}$ and $\textbf{H}$ magnetic
fields, which depend on the time $t$ and the space-variable $x$.
Moreover, we suppose that the initial electromagnetic field is
located into half-space $\mathbb{R}^N_{-} :=\{x=(x',x_N)\in
\mathbb{R}^N: x_N < 0\}$, i.\,e.
\begin{equation}\label{w}
\text{supp}\,w(0,.) \subset  \mathbb{R}^N_{-},
\end{equation}
where $w(x,t)= E^2 + H^2$.

Thus, the presented system $(M)$ is obtained from the classical
Maxwell's system ($M_0$) taking into account the state equations
(\ref{se0}) for isotropic nonlinear medium and Ohm's law for
current density (\ref{Ohm}). Mediums are describe to possess the
finite speed propagations property. There are many papers in which
energy decay was obtained for different problems concerning
Maxwell's equations. Well-posedness and asymptotic stability
results and decay of solutions are proved making use of different
techniques. Below, we mention some results concerning energy decay
and asymptotic of solutions.

Some linear evolution problems arise in the theory of hereditary
electromagnetism. Many authors studied the influence of
dissipation due to the memory on the asymptotic behavior of the
solutions (see \cite{B,FM, FM0,G,LN,LV,PR}). The polynomially
decay of the solutions when the memory kernel decays exponentially
or polynomially was shown in \cite{MNV}. It is studied the
asymptotic behavior of the solution of the linear problem
describing the evolution of the electromagnetic field inside a
rigid conducting material, whose constitutive equations contain
memory terms expressed by convolution integrals. These models were
proposed in \cite{NV} where it was shown that the exponential
decay of the memory kernel is able to produce a uniform rate decay
of energy in rigid conductors with electric memory.

The exact boundary controllability and stabilization of Maxwell's
equations have been studied by many authors (see \cite{NP} and
references therein). In \cite{NP} the internal stabilization of
Maxwell's equations with Ohm's law for space variable coefficients
is studied. Authors give sufficient conditions on parameters of
medium which guarantee the exponential decay of the energy of the
system. The result is based on observability estimate, obtained in
some particular cases by the multiplier method, a duality argument
and a weakening of norm argument, and argument used in internal
stabilization of scalar wave equations.

The energy decay of solutions of the scalar wave equation with
nonlinear damping in bounded domains has been shown in
\cite{D,H,N,S,TT,Z1,Z2}. In the case when there is no damping term
in the equation for the dielectric polarization, the long-time
asymptotic behavior of the solution of Maxwell's equations
involving generally nonlinear polarization and conductivity is
studied in \cite{J1}.

The propagation of electromagnetic waves in gas of quantum
mechanical system with two energy levels is considered in
\cite{J3}. The decay of the polarization field in a Maxwell--Bloch
system for $t\rightarrow \infty$ was shown.

The transient Landau-Lifschitz equations describing ferromagnetic
media without exchange interaction coupled with Maxwell's
equations is considered in \cite{J2}. The asymptotic behavior of
the solution of this mathematical model for micromagnetism is
studied. It is shown the strong convergence of the electromagnetic
field with respect to the energy norm for $t\rightarrow\infty$ on
bounded sets of nonvanishing electrical conductivity.

Following the dominant trend in the literature, we can conclude
that study of the system $(M)$ is not only of theoretical interest
but it is useful for applied researches. Since these authors are
not specialists in electromagnetism, we apologize in advance for
the omissions and inaccuracies. We hope that there is an
interdisciplinary audience which may find this useful, whether we
do not know any concrete mediums with proposed properties.

The present paper is organized as follows. In Section~2 we
formulate our main result. In Sections~3 we prove the finite speed
propagations property to some time, which depends on the
parameters of the problem and the initial electromagnetic field.
The method of proof is connected with nonhomogeneous variants of
Stampacchia lemma, in fact, it is an adaptation of local energy or
Saint--Venant principle like estimates method. Appendix~A contains
necessary interpolation inequalities and important properties of
nonhomogeneous functional inequalities.

\section{Main result}

We introduce the following concept of generalized solution of the
system $(M)$:

\begin{definition}
Let $n >1,\ p > 1,\ -p<m<p(n-1)$ and $w = E^2 + H^2$.  A pair
$(\textbf{E}(x,t), \textbf{H}(x,t))$ such that
$$
w \in C(0,T;L^1(\mathbb{R}^N)), \
 w^{\frac{m+p}{p}} \in
L^p(0,T;W^{1,p}(\mathbb{R}^N)), \ w_t \in L^1(Q_T)
$$
is called a solution to problem $(M)$ if for a.e. $t > 0$ the
integral identities
\begin{equation*}
\frac{1}{2} \int\limits_{\mathbb{R}^N}E^2(t,x)\eta(t,x)\,dx -
\frac{1}{2}\iint\limits_{Q_T} {E^2(t,x)\eta_t(t,x)\,dx\,dt}
+\iint\limits_{Q_T} {a_1 E^2(t,x)\eta(t,x) \,dx\,dt}
\end{equation*}
\begin{equation}\label{Def1}
 -
\iint\limits_{Q_T} {b_1\textbf{E}\, \mbox{curl}\, \textbf{H}
\,dx\,dt} = \frac{1}{2}
\int\limits_{\mathbb{R}^N}E^2(0,x)\eta(0,x)\,dx,
\end{equation}
\begin{equation*}
\frac{1}{2} \int\limits_{\mathbb{R}^N}H^2(t,x)\eta(t,x)\,dx-
\frac{1}{2}\iint\limits_{Q_T} {H^2(t,x)\eta_t(t,x)\,dx\,dt}
+\iint\limits_{Q_T} {a_2 H^2(t,x)\eta(t,x) \,dx\,dt}
\end{equation*}
\begin{equation}\label{Def2}
 + \iint\limits_{Q_T} {b_2 \textbf{H}\, \mbox{curl}\, \textbf{E}
\,dx\,dt} = \frac{1}{2}
\int\limits_{\mathbb{R}^N}H^2(0,x)\eta(0,x)\,dx,
\end{equation}
are satisfied for every $\eta \in C^1(Q_T)$.
\end{definition}

The main result is the following.

\begin{theorem}
Let the pair $(\textbf{E}(x,t), \textbf{H}(x,t))$ be a solution of
the problem $(M)$, in the sense of Definition~2.1. Let $p > 1,\ n
> 1$ $\bigl($and $n < 1 + \frac{(p - 1)(p + N)}{pN(2-p)}$ if $ p < 2\bigr)$, and
$$
\max \bigl\{- p,\, - p \bigl(1+ \tfrac{1}{N} -
\tfrac{n}{p}\bigr),\,
 - p\bigl(1+ \tfrac{1}{N} - \tfrac{n-1}{p-1}\bigr) \bigr\} < m < p (n - 2) + 1.
$$
Then there exists a time $ T^\ast > 0$, depending on known
parameters only (in particular, $\|w(x,0)\|_{L_1(\mathbb{R}^N)}$),
and a function $\Gamma(t) \in C[0, T] ,\ \Gamma(0) = 0$ such that
\begin{equation}\label{gg}
\Gamma(t) =  K\,\max \{ t^{\frac{p + N(m + p - n)}{p + N(m + p -
1)}},\ t^{\kappa}\} = K\,\left\{\begin{gathered} t^{\kappa}
\text{ for } t < 1,\\
t^{\frac{p + N(m + p - n)}{p + N(m + p - 1)}} \text{ for } t > 1,
\end{gathered}\right.
\end{equation}
where
$$
\kappa = \frac{p(p - 1 + N(m + p - n))[np + N(m + p - 1)]}{(p +
N(m + p - 1))[p(p(n-1) - m) + N(p-1)(m + p - 1)]},
$$
and
\begin{equation}\label{fsp}
\supp\,w(t,.) \subset \{x=(x',x_N)\in \mathbb{R}^N: x_N <
\Gamma(t)\} \ \forall\,0 < t < T^\ast,
\end{equation}
i.\,e. $\textbf{E}(x,t) = \textbf{H}(x,t) = 0$ for all $ x \in
\{x=(x',x_N)\in \mathbb{R}^N: x_N \geqslant \Gamma(t)\}$. Here $K
= K(n,m,p,N,\|w(0,x)\|_{L_1(\mathbb{R}^N)})$ is some positive
constant.
\end{theorem}

\begin{remark}
The statement of Theorem~1 stays true if we consider the problem
for system $(M)$ in some bounded domain. Then, instead of
(\ref{w}), we suppose that a support of initial energy of
electromagnetic field is contained in some ball into the domain.
\end{remark}

\section{Proof of finite speed of propagations}

Summing (\ref{Def1}) and (\ref{Def2}), in view of conditions
(\ref{c1}) and (\ref{c2}), we find that
\begin{multline}\label{I2}
\frac 12 \int\limits_{\mathbb{R}^N}w(t,x)\eta(t,x)\,dx -\frac 12
\iint\limits_{Q_T} {w(t,x)\eta_t(t,x)\,dx\,dt}
+ d_1\iint\limits_{Q_T} {w^{m} |\nabla\,w|^p \eta(t,x) \,dx\,dt} +\\
\iint\limits_{Q_T} {b_1 \,\text{div}\,(\textbf{E} \times
\textbf{H})\,\eta(x,t)\,dx\,dt}  \leqslant \frac 12
\int\limits_{\mathbb{R}^N}w(0,x)\eta(0,x)\,dx.
\end{multline}
Above we used the following relation:
\begin{equation}\label{R}
\text{div}\,(\textbf{E} \times \textbf{H}) = \textbf{H}
\,\text{curl}\,\textbf{E} - \textbf{E} \,\text{curl}\,\textbf{H}.
\end{equation}
From (\ref{I2}), (\ref{c3}) and (\ref{c4}) we get
\begin{multline}\label{I3-0}
\int\limits_{\mathbb{R}^N} {w(x,T)\, \eta(x,T) \,dx} -
\iint\limits_{Q_T} {w(x,t) \eta_t(x,t) \,dx\,dt} +
c \iint\limits_{Q_T} {|\nabla\,w^{\frac{m + p}{p}}|^p \eta(x,t) \,dx\,dt} \leqslant\\
\int\limits_{\mathbb{R}^N} {w(x,0)\, \eta(x,0) \,dx} + 2
d_2\iint\limits_{Q_T} {w^{n - 1}\,|\textbf{E} \times
\textbf{H}|\,|\nabla\,\eta(x,t)|\,dx\,dt} + \\
2 d_3 \iint\limits_{Q_T} {w^{n - 2}\,|\nabla\,w|\,|\textbf{E}
\times \textbf{H}| \eta(x,t)\,\,dx\,dt} \leqslant
\int\limits_{\mathbb{R}^N} {w(x,0)\, \eta(x,0) \,dx} + \\
\varepsilon \iint\limits_{Q_T} {|\nabla\,w^{\frac{m + p}{p}}|^p
\eta(x,t) \,dx\,dt} + c(\varepsilon) \iint\limits_{Q_T}
{w^{\frac{p(n - 1) - m}{p - 1}} \eta(x,t)\,dx\,dt} +\\
c\iint\limits_{Q_T} {w^n \,|\nabla\,\eta(x,t)|\,dx\,dt}
\end{multline}
for every nonnegative function $\eta(x,t) \in C^1 (Q_T)$, where
$\varepsilon > 0,\ p > 1,\ n > 1,\ - p < m < p(n -2) +1$ (i.\,e. $
\frac{p(n - 1) - m}{p - 1} >1$).

For an arbitrary $s \in \mathbb{R}^1$ and $\delta > 0$ we consider
the families of sets
$$
\begin{gathered}
\Omega (s) = \{x=(x',x_N)\in \mathbb{R}^N: x_N \geqslant s\},\
Q_T(s)= (0,T) \times \Omega (s) , \hfill\\
K (s,\delta ) = \Omega (s)\backslash \Omega (s + \delta ),\ K_T
(s,\delta ) = (0,T) \times  K (s,\delta ). \hfill
\end{gathered}
$$

Next we introduce our main cut-off functions $\eta _{s,\delta }
(x) \in C^1 (\mathbb{R}^N)$ such that $0 \leqslant \eta_{s,\delta
} (x) \leqslant 1 \  \forall  x \in \mathbb{R}^N$ and possess the
following properties:
$$
\eta _{s,\delta } (x) = \left\{
\begin{aligned} \hfill  0 \;
& , x \in \mathbb{R}^N \setminus \Omega(s),\\
\hfill  1\; & , x \in \Omega(s + \delta), \\
\end{aligned}\right. \ \ \  | \nabla\, \eta _{s,\delta }| \leqslant \tfrac{c} {\delta } \
\forall\,x \in K (s,\delta ).
$$
Choosing $\varepsilon > 0$ sufficiently small and
\begin{equation}\label{z}
\eta(x,t) = \eta _{s,\delta }(x) \exp \left( { - t \cdot T^{ - 1}
} \right)\ \forall\,T>0
\end{equation}
in integral inequality (\ref{I3-0}), we find
\begin{multline}\label{I3}
\mathop {\sup }\limits_{t \in (0,T)} \int\limits_{\Omega(s +
\delta)} {w(x,t)\,dx} + \frac{1}{T} \iint\limits_{Q_T (s +
\delta)} {w(x,t)\,dx\,dt} +
c \iint\limits_{Q_T(s + \delta)} {|\nabla\,w^{\frac{m + p}{p}}|^p \,dx\,dt} \leqslant\\
\int\limits_{\Omega (s)} {w(x,0)\,dx} + \frac{c}{\delta}
\iint\limits_{K_T (s,\delta)} {w^n \,dx\,dt} +
c\iint\limits_{Q_T(s)} {w^{\frac{p(n - 1) - m}{p - 1}}\,dx\,dt}=:
R_T(s,\delta),
\end{multline}
where $s \in \mathbb{R}^1 ,\ \delta > 0,\ T > 0$. Owing to
(\ref{w}), we have
\begin{equation}\label{w1}
\int\limits_{\Omega (s)} {w(x,0)\,dx} \equiv 0 \ \forall\,s
\geqslant 0.
\end{equation}
We introduce the functions related to $w(x,t)$:
$$
A_T (s): = \iint\limits_{Q_T (s)} {w^n\,dx\,dt},\ B_T (s): =
\iint\limits_{Q_T (s)} {w^{\frac{p(n - 1) - m}{p - 1}}\,dx\,dt}.
$$
Applying the interpolation inequality of Lemma~A.2 in the domain
$\Omega(s + \delta )$ to the function $v = w^{\frac{m + p}{p}}$
for $a = \frac{n\,p}{m + p},\ d = p,\ b = \frac{p}{m + p}$, $i =
0,\ j = 1$, and integrating the result with respect to time from
$0$ to $T$, we obtain
\begin{equation}\label{es1}
A_T(s + \delta) \leqslant c\, T^{1 - k_1} R_T^{1 +
\beta_1}(s,\delta),
\end{equation}
where  $k_1 = \tfrac{N(n - 1)}{p + N(m + p - 1)} < 1,\ \beta_1 =
\tfrac{p(n - 1)}{p + N(m + p - 1)}$, $m > n - p(1 + \frac{1}{N})$.
Similarly, applying the interpolation inequality of Lemma~A.2 in
the domain $\Omega(s + \delta )$ to the function $v = w^{\frac{m +
p}{p}}$ for $a = \frac{p(p(n-1) - m)}{(p-1)(m + p)},\ d = p,\ b =
\frac{p}{m + p}$, $i = 0,\ j = 1$, and integrating the result with
respect to time, we find that
\begin{equation}\label{es2}
B_T(s + \delta) \leqslant c\, T^{1 - k_2} R_T^{1 +
\beta_2}(s,\delta),
\end{equation}
where  $k_2 = \tfrac{N(p(n - 2) - m + 1)}{(p-1)(p + N(m + p - 1))}
< 1,\ \beta_2 = \tfrac{p(p(n - 2) - m + 1)}{(p-1)(p + N(m + p -
1))} $, $m > \frac{p(n-1)}{p -1} - p(1 + \frac{1}{N})$. Next we
define the function
$$
C_T(s) := (A_T (s))^{1 + \beta_2}+ (B_T (s))^{1 + \beta_1}.
$$
Then
\begin{equation}\label{D1}
C_T (s+\delta) \leqslant c\,F(T)\,\bigl[\delta^{- \beta}C_T^{1 +
\beta_1}(s) + C_T^{1 + \beta_2}(s) \bigr],
\end{equation}
where
$$
\beta = (1 + \beta_1)(1 + \beta_2),\ F(T)= \max \{ T^{(1 - k_1)(1
+ \beta_2)},\ T^{(1 - k_2)(1 + \beta_1)}\}.
$$
Below, we find some estimate $L^1$--norm of $w(x,t)$ by
$L^1$--norm of $w(x,0)$, which we will be used in the next
consideration.
\begin{lemma}
There exists some constant $c > 0$, depending on known parameters
of the problem, such that the following estimate
\begin{equation}\label{L1}
\int\limits_{\mathbb{R}^N} {w(x,t)\,dx}\leqslant
c\,\int\limits_{\mathbb{R}^N} {w(x,0)\,dx} \ \forall\, t \leqslant
T_1,
\end{equation}
is valid. Here $T_1$ depends on $m,\,p,\,n,\,N$ and
$\|w(x,0)\|_{L^1(\mathbb{R}^N)}$.
\end{lemma}

\begin{proof}
We set $s = - 2\delta ,\ \delta = s'
> 0$ in (\ref{I3}) and pass to the limit as $s' \to \infty $
\begin{multline}\label{I3-1}
\mathop {\sup }\limits_{t \in (0,T)} \int\limits_{\mathbb{R}^N}
{w(x,t)\,dx} + \frac{1}{T} \iint\limits_{Q_T} {w(x,t)\,dx\,dt} +
c \iint\limits_{Q_T} {|\nabla\,w^{\frac{m + p}{p}}|^p \,dx\,dt} \leqslant\\
\int\limits_{\mathbb{R}^N} {w(x,0)\,dx} + \iint\limits_{Q_T}
{w^{\frac{p(n - 1) - m}{p - 1}}\,dx\,dt}.
\end{multline}
Applying the interpolation inequality of Lemma~A.2 in
 $\mathbb{R}^N$ to the function $v = w^{\frac{m + p}{p}}$ for $a =
\frac{p(p(n-1) - m)}{(m + p)(p-1)},\ d = p,\ b = \frac{p}{m + p}$,
$i = 0,\ j = 1$, and Young's inequality, we find that
\begin{multline}\nonumber
\int\limits_{\mathbb{R}^N} {w^{\frac{p(n - 1) - m}{p - 1}}\,dx}
\leqslant c \left ( \int\limits_{\mathbb{R}^N}
{|\nabla\,w^{\frac{m+p}{p}}|^p\,dx}\right)^\frac{a\theta}{p}\left(
\int\limits_{\mathbb{R}^N} {w}\,dx\right)^{\frac{a(1-\theta)}{b}}
\leqslant\\
 \varepsilon \int\limits_{\mathbb{R}^N}
{|\nabla\,w^{\frac{m+p}{p}}|^p\,dx}+c(\varepsilon)\left(
\int\limits_{\mathbb{R}^N}
{w\,dx}\right)^{\frac{ap(1-\theta)}{b(p-a\theta)}}  \
\forall\,\varepsilon > 0,
\end{multline}
where $ \theta=\frac{N(m+n)(p(n-2)-m+1)}{(N(m+p-1)+p)(p(n-1)-m)}$.
Integrating this inequality with respect to time from $0$ to $T$,
we obtain
\begin{multline}\label{I3-2}
\iint\limits_{Q_T} {w^{\frac{p(n - 1) - m}{p - 1}}\,dx\,dt}
\leqslant  \varepsilon \iint\limits_{Q_T}
{|\nabla\,w^{\frac{m+p}{p}}|^p\,dx}+c(\varepsilon)\int\limits_0^T
{\left( \int\limits_{\mathbb{R}^N}
{w\,dx}\right)^{\frac{ap(1-\theta)}{b(p-a\theta)}}dt}.\\
\end{multline}
Choosing $\varepsilon > 0$ sufficiently small, from
\eqref{I3-1},\eqref{I3-2} we have
\begin{multline}\label{I3-3}
\mathop {\sup }\limits_{t \in (0,T)} \int\limits_{\mathbb{R}^N}
{w(x,t)\,dx} + \frac{1}{T} \iint\limits_{Q_T} {w(x,t)\,dx\,dt} +
c \iint\limits_{Q_T} {|\nabla\,w^{\frac{m + p}{p}}|^p \,dx\,dt} \leqslant\\
\int\limits_{\mathbb{R}^N} {w(x,0)\,dx} +
c\,\int\limits_0^T{\left( \int\limits_{\mathbb{R}^N}
w\,dx\right)^{\frac{ap(1-\theta)}{b(p-a\theta)}}dt}.
\end{multline}
From the last inequality we deduce that for every $t: 0<t<T$ the
following inequality is valid
$$
\int\limits_{\mathbb{R}^N} {w(x,t)\,dx}\leqslant
\int\limits_{\mathbb{R}^N} {w(x,0)\,dx} +c\,\int\limits_0^t{\left(
\int\limits_{\mathbb{R}^N} w(x,\tau)\,dx\right)^{\gamma}\!\!\!
d\tau},
$$
where $ \gamma =\frac{(N-1)(p(n-1)-m) +N(p-1)(m+p)} {p( p -1 +
N(m+p-n))}$. Applying Lemma~A.3 from Appendix~A we obtain
\eqref{L1} with
\begin{equation}\label{T1}
T_1 := \left\{\begin{gathered}
 \tfrac{2}{1 - \gamma} \biggl( \int\limits_{\mathbb{R}^N}
w(x,0)\,dx \biggr)^{1-\gamma} \text{ if } \gamma < 1, \hfill\\
\tfrac{1}{2(\gamma -1)} \biggl( \int\limits_{\mathbb{R}^N}
w(x,0)\,dx \biggr)^{\gamma - 1} \text{ if } \gamma > 1, \hfill\\
\end{gathered}\right.
\end{equation}
and $T_1 \to 0$ as $\|w(x,0)\|_{L_1(\mathbb{R}^N)}\rightarrow 0$.
\end{proof}

Further, using the definition of the functions $C_T (s)$ and
\eqref{L1}, we get
\begin{equation}\label{J0}
C_T (s_0) \leqslant K_0\,F(T) \ \forall\, T \leqslant T_1.
\end{equation}
where the positive constant $K_0$ depends on $n,\ m,\ p,\ N$ and $
\|w(x,0)\|_{L^1(\mathbb{R}^N)}$.

Now we choose the parameter $\delta > 0$ which was arbitrary up to
now:
$$
\delta _T (s): = \left[ {\frac{2c}{1 - H_T(s_0)}\,F(T)
\,C_T^{\beta_1} (s)} \right]^{\frac{1} {\beta}},
$$
where the function $ H_T(s) = c\,F(T)\,C_T^{\beta_2}(s)$ is such
that $H_T(s_0) < 1$ at some point $s_0 \geqslant 0$, whence we get
that
\begin{equation}\label{T2}
T \leqslant T_2 = c \min
\{K_0^{-\frac{\beta_2}{(1-k_1)(1+\beta_2)^2}},
K_0^{-\frac{\beta_2}{(1-k_2)(1+\beta_1)(1+\beta_2)}}\},
\end{equation}
and $ T_2 \to \infty$ as
$\|w(x,0)\|_{L_1(\mathbb{R}^N)}\rightarrow 0$.

We obtain the following main functional relation for the function
$\delta_T (s)$:
\begin{equation}\label{J}
\delta_T (s +\delta_T (s)) \leqslant \varepsilon \, \delta_T (s)\
\forall \, s \geqslant s_0 \geqslant 0, \ 0 < \varepsilon =
\bigl(\tfrac{1 + H_T(s_0)}{2}\bigr)^{\frac{\beta_1} {\beta}} < 1
\end{equation}
$\forall\,\,\,0 < T < T^\ast := \min \{T_1, T_2\}$, where $T_1$ of
\eqref{T1} and $T_2$ of \eqref{T2}. Now we apply Lemma~A.1 to the
function $\delta_T (s)$ of (\ref{J}). As a result, we obtain
\begin{equation}\label{J1}
\delta_T (s) \equiv 0 \ \forall\,s \geqslant s_ 0 + \tfrac{1}{1 -
\varepsilon} \delta_T(s_0).
\end{equation}
Then, in view of (\ref{J0}), we find
\begin{multline}\nonumber
\delta_T (s_0)\leqslant c\,
[C_T^{\beta_1}(s_0)F(T)]^{\frac{1}{\beta}} \leqslant c\,
[F^{1+\beta_1}(T)]^{\frac{1}{\beta}} \leqslant c
\,(F(T))^{\frac{1}{1+\beta_2}}= \\
c \,\max \{ T^{1 - k_1},\
T^{\frac{(1 - k_2)(1 + \beta_1)}{1 + \beta_2}}\}
\end{multline}
$\forall\,0 < T < T^\ast$. Choosing in (\ref{J1}) $s_0 = 0$ and
$$
s = \Gamma (T) = c\,\,\max \{ T^{\frac{p + N(m + p - n)}{p + N(m +
p - 1)}},\ T^{\kappa}\} = c\, \left\{\begin{gathered} T^{\kappa}
\text{ for } T < 1,\\
T^{\frac{p + N(m + p - n)}{p + N(m + p - 1)}} \text{ for } T > 1
\end{gathered}\right.
$$
$\forall\,0 < T < T^\ast$, $ \kappa = \frac{p(p - 1 + N(m + p -
n))[np + N(m + p - 1)]}{(p + N(m + p - 1))[p(p(n-1) - m) +
N(p-1)(m + p - 1)]}$. Thus $w(T,x) \equiv 0$ for all $x \in
\{x=(x',x_N)\in \mathbb{R}^N: x_N \geqslant \Gamma(t)\}$. And
Theorem~1 is proved completely. $\square$

\section*{Appendix A}
\renewcommand{\thesection}{A}\setcounter{equation}{0}
\setcounter{lemma}{0}
\renewcommand{\thetheorem}{A.\arabic{lemma}}

\begin{lemma}$\cite{SH2}$
Let the nonnegative continuous nonincreasing function $f(s):
[s_0,\infty) \to \mathbb{R}^1$ satisfies the following functional
relation:
$$
f(s + f(s)) \leqslant \varepsilon\,f(s)\ \forall\, s \geqslant
s_0,\ 0 < \varepsilon < 1.
$$
Then $f(s)\equiv 0 \ \forall\, s \geqslant s_0+(1 -
\varepsilon)^{-1}f(s_0)$.
\end{lemma}

\begin{lemma}$\cite{Ni}$
If $\Omega  \subset \mathbb{R}^N $ is a bounded domain with
piecewise-smooth boundary, $a > 1$, $b \in (0, a),\ d
> 1,$ and $0 \leqslant i < j,\ i,j \in \mathbb{N}$, then there exist positive constants
$d_1$ and $d_2$ $(d_2 = 0$ if the domain $\Omega$ is unbounded$)$
that depend only on $\Omega ,\ d,\ j,\ b,$ and $N$ and are such
that, for any function $v(x) \in W_d^j (\Omega ) \cap L^b (\Omega
)$, the following inequality is true:
$$
\left\| {D^i v} \right\|_{L^a (\Omega )}  \leqslant d_1 \left\|
{D^j v} \right\|_{L^d (\Omega )}^\theta  \left\| v \right\|_{L^b
(\Omega )}^{1 - \theta }  + d_2 \left\| v \right\|_{L^b (\Omega )}
$$
where $\theta  = \frac{{\tfrac{1} {b} + \tfrac{i} {N} - \tfrac{1}
{a}}} {{\tfrac{1} {b} + \tfrac{j} {N} - \tfrac{1} {d}}} \in \left[
{\tfrac{i} {j},1} \right)$.
\end{lemma}

\begin{lemma}$\cite{Bi}$
Suppose that $v(t)$ is a nonnegative summable function on $[0,T]$
that, for almost all $t \in[0,T],$ satisfies the integral
inequality
$$
 v(t) \leqslant k + m
\int\limits_0^t {h(\tau) g(v(\tau))\,d\tau }
$$
 where $k \geqslant
0, m \geqslant 0$, $h(\tau)$ is summable on $[0,T],$ and $g(\tau)$
is a positive function for $\tau
> 0$. Then
$$
 v(t) \leqslant G^{-1}\Biggl( {G(k) + m \int\limits_0^t
{h(\tau)\,d\tau }} \Biggr)
$$
 for almost all $t \in[0,T]$. Here
$G(v)= \int \limits_{v_0}^{v} {\frac{d\tau}{g(\tau)}},\ v > v_0 >
0$.
\end{lemma}

\end{document}